\documentclass{jT}

\usepackage{amssymb,amsmath,latexsym,amsthm}
\usepackage[english]{babel}

\theoremstyle{definition}

\renewcommand{\le}{\leqslant}
\renewcommand{\ge}{\geqslant}
\numberwithin{equation}{section}

\def\DJ{\leavevmode\setbox0=\hbox{D}\kern0pt\rlap
{\kern.04em\raise.188\ht0\hbox{-}}D}

\def\hf{{\textstyle\frac{1}{2}}}
\def\a{\alpha}\def\b{\beta}

\def\e{\varepsilon}
\def\f{\varphi}
\def\G{\Gamma}

\def\s{\sigma}
\def\t{\theta}
\def\={\;=\;}

\def\zt{\zeta(\hf+it)}

\def\D{\Delta}

\def\R{\Re{\rm e}\,} 
\def\z{\zeta}

\def\f{\varphi}

\font\tenmsb=msbm10
\font\sevenmsb=msbm7
\font\fivemsb=msbm5
\newfam\msbfam
\textfont\msbfam=\tenmsb
\scriptfont\msbfam=\sevenmsb
\scriptscriptfont\msbfam=\fivemsb
\def\Bbb#1{{\fam\msbfam #1}}

\def \NN {\Bbb N}

\font\teneufm=eufm10
\font\seveneufm=eufm7
\font\fiveeufm=eufm5
\newfam\eufmfam
\textfont\eufmfam=\teneufm
\scriptfont\eufmfam=\seveneufm
\scriptscriptfont\eufmfam=\fiveeufm
\def\mathfrak#1{{\fam\eufmfam\relax#1}}

\def\Residuo{\mathop{\hbox{\rm {Res}}}}

\def \NN {\Bbb N}

\def\D{\Delta}
\def\a{\alpha}
\def\b{\beta} \def\e{\varepsilon}

\def\d{\,{\rm d}}

\begin{document}

\title[Mean square of the divisor function in short intervals]
{On the mean square of the divisor function in short intervals}


\author[Aleksandar {\sc ivi\'c}]{{\sc Aleksandar} IVI\'C}
\address{Aleksandar {\sc ivi\'c}\\
Katedra Matematike RGF-a\\ Universitet u Beogradu,  \DJ u\v sina 7\\
11000 Beograd, Serbia}
\email{ivic@rgf.bg.ac.yu, aivic\_2000@yahoo.com}

\maketitle

\begin{resume}
On donne des estimations pour la moyenne quadratique de
$$
\int_X^{2X}\left(\D_k(x+h) - \D_k(x)\right)^2\d x,
$$
o\`u $h = h(X)\gg1,\; h = o(x)\;
{\rm{quand}}\;X\to\infty$ et
$h$ se trouve dans un intervalle convenable. Pour $k\ge2$ un entier fix\'e,
$\D_k(x)$ et le terme d'erreur pour la fonction
sommatoire  de la fonction des  diviseurs  $d_k(n)$,
gener\'ee par $\z^k(s)$.
\end{resume}

\begin{abstr}
We provide upper bounds for the mean square integral
$$
\int_X^{2X}\left(\D_k(x+h) - \D_k(x)\right)^2\d x,
$$
where $h = h(X)\gg1,\; h = o(x)\;
{\rm{as}}\;X\to\infty$ and
$h$ lies in a suitable range. For $k\ge2$ a fixed integer,
$\D_k(x)$ is the error term in the asymptotic formula
for the summatory function of the divisor function $d_k(n)$,
generated by $\z^k(s)$.
\end{abstr}

\bigskip
\section{Introduction}

Let, for a fixed integer $k\ge2$, $d_k(n)$ denote the (generalized)
divisor function which denotes the number of ways $n$ can be written
as a product of $k$ factors. This is a well-known multiplicative
function ($d_k(mn) = d_k(m)d_k(n)$ for coprime $m,n\in\NN$, and $d_2(n)\equiv d(n)$
is the classical number of divisors function). Besides this definition one
has the property that $d_k(n)$ is generated by $\z^k(s)$, where $\z(s)$
is the Riemann zeta-function, defined as
$$
\z(s) = \sum_{n=1}^\infty n^{-s}\qquad(\s = \R s > 1),
$$
and otherwise by analytic continuation. Namely
$$
\z^k(s) = {\left(\sum_{n=1}^\infty n^{-s}\right)}^k
= \sum_{n=1}^\infty d_k(n)n^{-s}\qquad(\s = \R s > 1),
$$
and this connects the problems related to $d_k(n)$ to zeta-function
theory. Let as usual
$$
\sum_{n\le x}d_k(n) = \Residuo_{s=1}\,{x^s\z^k(s)\over s} + \D_k(x)
= xP_{k-1}(\log x) + \D_k(x),
\leqno(1.1)
$$
where $P_{k-1}(t)$ is a polynomial of degree $k-1$ in $t$,
all of whose coefficients can be calculated explicitly
(e.g., $P_1(t) = \log t + 2\gamma-1$, where $\gamma = -\G'(1)= 0.577215\ldots$
is Euler's constant). Thus $\D_k(x)$ represents the error term
in the asymptotic formula for the summatory function of $d_k(n)$,
and a vast literature on this subject exists (see e.g., [2] or [6]).
Here we shall be concerned with the ``short difference"
$$
\D_k(x+h) - \D_k(x)\qquad(1 \ll h \ll x,\; h=o(x)\; {\rm{as}}\; x\to\infty).
\leqno(1.2)
$$
 The meaning of
``short" comes from the condition $h = o(x)$ as $x\to\infty$, so that $h$ is indeed much
smaller in comparison with $x$.

\medskip
As in analytic number theory one is usually interested in the
averages of error terms, where the averaging usually smoothens
the irregularities of distribution of the function in question, we
shall be interested in mean square estimates of (1.2), both discrete and
continuous. To this end we therefore define
the following means (to stress the analogy between the discrete
and the continuous, $x$ is being kept
as both the continuous and the integer variable):
$$
{\sum}_k(X,h) :\;=\; \sum_{X\le x\le2X}{\Bigl(\D_k(x+h) - \D_k(x)\Bigr)}^2,
\leqno(1.3)
$$
$$
I_k(X,h) :\;=\; \int_X^{2X}{\Bigl(\D_k(x+h) - \D_k(x)\Bigr)}^2\d x.
\leqno(1.4)
$$
The problem is then to find non-trivial upper
bounds for (1.3)--(1.4), and to show
for which ranges of $h = h(k,X)\;(=o(X)\,)$ they are valid.
Theorem 1 (see Section 2) provides some results in this direction.
Namely the ``trivial"
bound in all cases is the function $X^{1+\e}h^2$, where here and later $\e>0$
denotes arbitrarily small positive constants, not necessarily the same ones
at each occurrence. This comes from the elementary
bound $d_k(n) \ll_\e n^\e$ and ($1\ll h \ll x$)
$$
\leqno(1.5)
$$
\begin{eqnarray*}
\D_k(x+h) - \D_k(x)&=& \sum_{x<n\le x+h}d_k(n) + xP_{k-1}(\log x)\\&
-& (x+h)P_{k-1}(\log(x+h))\\&
\ll_\e& x^\e\sum_{x<n\le x+h}1 + h\log^{k-1}x \ll_\e x^\e h,
\end{eqnarray*}
where we used the mean value theorem.
Here and later $\e\,(>0)$ denotes arbitrarily small constants,
not necessarily the same ones at each occurrence, while $a \ll_\e b$
(same as $a = O_\e(b)$) means that the implied constant depends on $\e$.
Note that, from the work of
P. Shiu [5] on multiplicative functions, one has the bound
$$
\sum_{x<n\le x+h}d_k(n) \ll h\log^{k-1}x\qquad(x^\e \le h \le x),
$$
hence in this range (1.5) can be improved a bit, and thus we can also consider
$$
\D_k(x+h) - \D_k(x) \ll h\log^{k-1}x\qquad(x^\e \le h \le x)\leqno(1.6)
$$
as the ``trivial"  bound. It should be mentioned that the cases $k=2$ of
(1.3) and (1.4) have been treated by Coppola--Salerno [1] and M. Jutila [4],
respectively, so that we shall concentrate here on the case when $k>2$. In these
papers it had been shown that, for $X^\e \le h \le \hf \sqrt{X},\, L = \log X$,
$$
{\sum}_2(X,h) = \frac{8}{\pi^2}Xh\log^3\Bigl(\frac{\sqrt{X}}{ h}\Bigr)
+ O(XhL^{5/2}\sqrt{L}),\leqno(1.7)
$$
$$
\leqno(1.8)
$$
\begin{eqnarray*}
I_2(X,h) &=&
{1\over4\pi^2}\sum_{n\le {X\over2h}}{d^2(n)\over n^{3/2}}\int\limits_X^{2X}
x^{1/2}\left|\exp\left(2\pi ih\sqrt{{n\over x}}\,\right)-1\right|^2\d x\\&
+ &O_\e(X^{1+\e}h^{1/2}).
\end{eqnarray*}
From (1.8) Jutila deduces ($a\asymp b$ means that $a \ll b \ll a$) that
$$
I_2(X,h) \;\asymp\;Xh\log^3\Bigl({\sqrt{X}\over h}\Bigr)
\qquad(X^\e \le h \le  X^{1/2-\e}),
\leqno(1.9)
$$
but it is not obvious that $I_2(X,h)$ is asymptotic to the main term on the
right-hand side of (1.7). This, however, is certainly true, and will follow
from our Theorem 2 (see Section 2) and from (1.7).
Theorem 2 says that, essentially, the sums
${\sum}_k(X,h)$ and $I_k(X,h)$ are of the same order of magnitude.
It is also true that, for $X^\e \le h \le \hf \sqrt{X},\, L = \log X$,
$$
\int_X^{2X}\Bigl(E(x+h) - E(x)\Bigr)^2\d x =
{8\over\pi^2}Xh\log^3\Bigl({\sqrt{X}\over h}\Bigr)
+ O(XhL^{5/2}\sqrt{L}),\leqno(1.10)
$$
implying in particular that
$$
E(x+h) - E(x) = \Omega \Biggl\{\sqrt{h}\log^{3/2}\Bigl({\sqrt{x}\over h}\Bigr)\Biggr\}
\qquad(x^\e \le h \le  x^{1/2-\e}).\leqno(1.11)
$$
Here, as usual,
$$
E(T) := \int_0^T|\zt|^2\d t - T\left(\log{T\over2\pi} + 2\gamma -1\right)
$$
represents the error term in the mean square formula for $|\zt|$ (see e.g.,
Chapter 15 of [2] for a comprehensive account), while $f(x) = \Omega(g(x))$
means that $\lim_{x\to\infty}f(x)/g(x)\ne0$.
 Namely Jutila (op. cit.) has shown that the integral in (1.10) equals
the expression on the right-hand side of (1.8), hence the conclusion follows from
 Theorem 2 and the above discussion. The omega-result (1.11) shows that
the difference $E(x+h) - E(x)$ cannot be too small in a fairly wide range for $h$.
An omega-result analogous to (1.11) holds for $\D(x+h)-\D(x)$ as well, namely
$$
\D(x+h) - \D(x) = \Omega \Biggl\{\sqrt{h}\log^{3/2}\Bigl({\sqrt{x}\over h}\Bigr)\Biggr\}
\qquad(x^\e \le h \le  x^{1/2-\e}).
\leqno(1.12)
$$

\medskip
Concerning the true order of $\D_k(x+h) - \D_k(x)$, we remark that on the basis
of (1.9) M. Jutila [4] conjectured that
$$
\D(x+h) - \D(x) \ll_\e x^\e\sqrt{h}\qquad(x^\e \le h \le x^{1/2-\e}),\leqno(1.13)
$$
which would be close to best possible, in view of (1.12). The range
$x^\e \le h \le x^{1/2-\e}$ is essential here, since for $h$ much larger than
$x^{1/2}$, one expects $\D(x+h)$ and $\D(x)$ to behave like independent random variables,
and in that case the quantities in question may not be ``close" to one another.
Perhaps one has (1.13) for $\D(x)$ replaced by $\D_k(x)$ in a suitable range of $h$
as well, but this is a difficult question.

\medskip
Further sharpenings of (1.7), (1.9) and (1.10) were recently obtained by the author in [3].
Namely, for $1 \ll U = U(T) \le \hf {\sqrt{T}}$ we have ($c_3 = 8\pi^{-2}$)
$$
\leqno(1.14)
$$
\begin{eqnarray*}
\int_T^{2T}\Bigl(\D(x+U)-\D(x)\Bigr)^2\d x & =& TU\sum_{j=0}^3c_j\log^j
\Bigl({\sqrt{T}\over U}\Bigr) \\&
+& O_\e(T^{1/2+\e}U^2) + O_\e(T^{1+\e}U^{1/2}),
\end{eqnarray*}
and the result remains true if $\D(x+U)-\D(x)$ is replaced by $E(x+U)-E(x)$
(with different $c_j$'s).
From (1.14) and Theorem 2 one obtains that a formula analogous to (1.14)
holds also for $\sum_2(X,h)$ in (1.7).

\medskip
\section{Statement of results}
\medskip
First we define $\s(k)$ as a number satisfying $\hf\le\s(k)<1$, for which
$$
\int_0^T|\z(\s(k)+it)|^{2k}\d t \;\ll_\e\;T^{1+\e}\leqno(2.1)
$$
holds for a fixed integer $k\ge2$.
From zeta-function theory (see [2], and in particular Section 7.9  of [6]) it is
known that such a number exists for any given $k\in\NN$, but it is  not uniquely
defined, as one has
$$
\int_0^T|\z(\s+it)|^{2k}\d t \;\ll_\e\;T^{1+\e}\qquad(\s(k)\le\s<1).\leqno(2.2)
$$
From Chapter 8 of [2] it follows that
 one has $\s(2) = \hf, \s(3) = {7\over12}, \s(4) = {5\over8},
\s(5) = {41\over60}$ etc., but it is not easy to write down (the best known
value of) $\s(k)$ explicitly as a function of $k$. Note that the famous,
hitherto unproved Lindel\"of hypothesis that $\zt \ll_\e |t|^\e$ is equivalent
to the fact that $\s(k) = \hf\; (\forall k \in \NN)$. Our aim is to
find an upper bound for $I_k(X,h)$ in (1.4) which is better than the
the trivial bound $O_\e(X^{1+\e}h^2)$ by the use of (1.5), or
$O(Xh^2\log^{2k-2}X)$, by the use of (1.6). Now we can formulate

\bigskip
THEOREM 1. {\it Let $k\ge3$ be a fixed integer. If $\s(k) = \hf$ then,
for $X^\e \le h = h(X) \le X^{1-\e}$,
$$
\int_X^{2X}{\Bigl(\D_k(x+h)-\D_k(x)\Bigr)}^2\d x \ll_\e X^{1+\e}h^{4/3}.
\leqno(2.3)
$$
If $\hf < \s(k) < 1$, and $\t(k)$ is any constant satisfying $2\s(k)-1< \t(k)<1$,
then there exists $\e_1 = \e_1(k)>0$ such that,
for $X^{\t(k)} \le h = h(X) \le X^{1-\e}$,}
$$
\int_X^{2X}{\Bigl(\D_k(x+h)-\D_k(x)\Bigr)}^2\d x \ll_{\e_1} X^{1-\e_1}h^{2}.
\leqno(2.4)
$$

\bigskip
{\bf Corollary 1}. If the Lindel\"of hypothesis is true, then
the bound (2.3) holds
for all $k$.

\medskip
{\bf Corollary 2}. From the known values of $\s(k)$ mentioned above it transpires
that one may unconditionally take $\t(3) = {1\over6}+\e,
\t(4) = {1\over4}+\e, \t(5) = {11\over30}+\e$, etc.

\medskip
{\bf Remark 1}. The bound in (2.3) is fairly sharp, while the bound in (2.4)
is a little better than the trivial bound $Xh^2(\log X)^{2k-2}$ (cf. (1.6)).

\medskip
{\bf Remark 2}. Theorem 1 holds also for $k=2$, but in this case
a sharper result follows from (1.9) in the range $X^\e \le h \le X^{1/2-\e}$.
It could be true that, for $k>2$ fixed,  the weak analogue of (1.7),
namely the bound $X^{1+\e}h$ holds in a suitable range for $h$ (depending on $k$),
but this seems unattainable at present.

\medskip
{\bf Remark 3}. Note that  the integrals in (2.3)--(2.4)
are trivially bounded by $X^{1+2\b_k+\e}$, where as usual
$$
\b_k \;:=\; \inf\Bigl\{\;b_k\;:\; \int_1^X\D^2_k(x)\d x \ll X^{1+2b_k}\;\Bigr\}
$$
for fixed $k\ge2$. It is known (see Chapter 13 of [2])
that $\b_k \ge (k-1)/(2k)$ for every $k$,
$\b_k = (k-1)/(2k)$ for $k = 2,3,4$, $\b_5 \le 9/20$ (see W. Zhang [7]),
$\b_6\le \hf$, etc. This gives an insight when Theorem 1 gives a non-trivial result.

\medskip

Our second result is primarily a technical one.
It establishes the connection between the discrete
means ${\sum}_k(X,h)$ (see (1.3)) and its
continuous counterpart $I_k(X,h)$ (see (1.4)),
precisely in the range where we expect the $\D$-functions to be
close to one another. This is

\bigskip
THEOREM 2. {\it For $1 \ll h = h(X) \le \hf\sqrt{X}$ we have
$$
{\sum}_2(X,h) = I_2(X,h) + O(h^{5/2}\log^{5/2}X),\leqno(2.5)
$$
while, for a fixed integer $k\ge3$,}
$$
{\sum}_k (X,h) = I_k(X,h) + O_\e(X^\e h^{3}).\leqno(2.6)
$$

\bigskip
\section{Proof of Theorem 1}
\medskip
We start from Perron's classical inversion formula
(see e.g., the Appendix of [2]). Since $d_k(n) \ll_\e n^\e$, this yields
$$
\frac{1}{2\pi i}\sum_{n\le x}d_k(n) =
\int_{1+\e-iT}^{1+\e+iT}\frac{x^s}{s}\,\z^k(s)\d s
+ O_\e(X^{1+\e}T^{-1}),\leqno(3.1)
$$
where $X\le x \le 2X$, and $T$ is parameter satisfying $1\ll T \ll X$ that
will be suitably chosen a little later.
We replace the segment of integration by the contour joining the points
$$
1+\e-iT,\, \s(k)-iT,\,\s(k)+iT,\,1+\e+iT.
$$
In doing this we encounter the pole of $\z^k(s)$ at $s=1$ of order $k$, and the
residue at this point will furnish $xP_{k-1}(\log x)$, the main term in (1.1).
Hence by the residue theorem (3.1) gives,
applied once with $x$ and once with $x+h$,
\pagebreak
$$
\leqno(3.2)
$$
\begin{eqnarray*}
\D_k(x+h) - \D_k(x)&=& {1\over2\pi i}\int_{\s(k)-iT}^{\s(k)+iT}
{(x+h)^s-x^s\over s}\,
\z^k(s)\d s \\&
\,+ &O_\e(X^{1+\e}T^{-1}) + O(R_k(x,T)),
\end{eqnarray*}
where
$$
R_k(x,T) := {1\over T}\int_{\s(k)}^{1+\e}x^\a|\z(\a+iT)|^k\d \a.\leqno(3.3)
$$
By using the Cauchy-Schwarz inequality for integrals, (2.1) and (2.2), it follows that
$$
\leqno(3.4)
$$
\begin{eqnarray*}
\int_{T_0}^{2T_0}R_k(x,T)\d T &
\ll& {1\over T_0}\int_{\s(k)}^{1+\e}x^\a\left(\int_{T_0}^{2T_0}
|\z(\a+iT)|^k\d T\right)\d \a\\&
\ll_\e& {X^{1+\e}\over T_0}\sup_{\s(k)\le\a\le1+\e}
\left(T_0\int_{T_0}^{2T_0}|\z(\a+iT)|^{2k}\d T\right)^{1/2}\\&
\ll_\e& X^{1+\e}.
\end{eqnarray*}
Therefore (3.4) implies that there exists $T\in [T_0,\,2T_0]$ such that
$$
R_k(x,T) \,< \,cX^{1+\e}T_0^{-1}\leqno(3.5)
$$
for a suitable $c>0$, uniformly in $X\le x \le 2X$. It is this $T$ that we initially
take in (3.2)--(3.3), and using
$$
{(x+h)^s-x^s\over s} = \int_0^h (x+v)^{s-1}\d v\qquad(s\ne 0)
$$
we obtain from (3.2) and (3.5)
$$
\D_k(x+h) - \D_k(x)= {1\over2\pi i}\int_{\s(k)-iT}^{\s(k)+iT}\int_0^h (x+v)^{s-1}\d v
\,\z^k(s)\d s + O_\e(X^{1+\e}T_0^{-1}).\leqno(3.6)
$$
On squaring (3.6) and integrating over $x$, we obtain
$$
\leqno(3.7)
$$
\begin{eqnarray*}
&{}&
\int_X^{2X}{\left(\D_k(x+h)-\D_k(x)\right)}^2\d x\\&
\ll_\e& \int_X^{2X}{\Bigl|\int_{-T}^T\int_0^h (x+v)^{\s(k)-1+it}
\z^k(\s(k)+it)\d v\d t\Bigl|}^2\d x
+ X^{3+\e}T_0^{-2}.
\end{eqnarray*}
Let now $\f(x)\;(\ge0)$ be a smooth function supported in $[X/2,\,5X/2]$, such that
$\f(x) = 1$ when $X \le x \le 2X$ and $\f^{(r)}(x) \ll_r X^{-r}\;(r=0,1,2,\ldots\,)$.
In the integrals under the absolute value signs in (3.7) we exchange the order
of integration and then
use the Cauchy-Schwarz inequality for integrals. We infer that the integral
on the right-hand side of (3.7) does not exceed
\begin{eqnarray*}
&{}&
h\int_{X/2}^{5X/2}\f(x)\int_0^h{\Bigl|\int_{-T}^T (x+v)^{\s(k)-1+it}
\z^k(\s(k)+it)\d t\Bigl|}^2\d v\,\d x\\&
=& h\int_0^h\int_{-T}^T\int_{-T}^T\z^k(\s(k)+it)\z^k(\s(k)-iy)J\,\d y\,\d t\,\d v,
\end{eqnarray*}
say, where
$$
J = J_k(X;v,t,y) :=\int_{X/2}^{5X/2}\f(x)(x+v)^{2\s(k)-2}(x+v)^{i(t-y)}\d x.
$$
Integrating by parts we obtain, since $\f(X/2) = \f(5X/2) = 0$,
$$
J = {-1\over i(t-y)+1}\int_{X/2}^{5X/2}(x+v)^{2\s(k)-1+i(t-y)}\Bigl(\f'(x)
+ {2\s(k)-2\over x+v}\f(x)\Bigr)\d x.
$$
By repeating this process it is seen that each time our integrand will be
decreased by the factor of order
$$
\ll \; {X\over|t-y|+1}\cdot {1\over X} \;\ll_\e\; X^{-\e}
$$
for $|t-y| \ge X^\e$. Thus if we fix any $\e, A>0$, the contribution of
$|t-y|\ge X^\e$ will be $\ll X^{-A}$ if we integrate by parts $r = r(\e,A)$
times. For $|t-y|\le X^\e$ we estimate the corresponding contribution to $J$
trivially to obtain that the integral on the right-hand side of (3.7) is
\begin{eqnarray*}
&\ll_\e& h^2X^{2\s(k)-1}\int_{-T}^T\int_{-T,|t-y|\le X^\e}^T
|\z(\s(k)+it)\z(\s(k)+iy)|^k
\d y\,\d t\\&
\ll_\e& h^2X^{2\s(k)-1}\int_{-T}^T|\z(\s(k)+it)|^{2k}
\left(\int_{t-X^\e}^{t+X^\e}\d y\right)\d t
\\&
\ll_\e& h^2X^{2\s(k)-1+\e}T,
\end{eqnarray*}
where we used (2.1) and the elementary inequality
$$
|ab| \le \hf(|a|^2+|b|^2).
$$
Since $T_0 \le T\le 2T_0$, it is seen that the left-hand side of (3.7) is
$$
\ll_\e\;X^\e(h^2T_0X^{2\s(k)-1} + X^3T_0^{-2}).\leqno(3.8)
$$

\medskip
If $\s(k) = \hf$, then
$$
h^2T_0 + X^3T_0^{-2} = 2Xh^{4/3}
$$
with the choice $T_0 = Xh^{-2/3}$, which clearly satisfies $1 \ll T_0 \ll X$. Therefore
(2.3) follows from (3.7) and (3.8).

\medskip
If $\hf < \s(k) <1$, then we choose first
$$
T_0 \= X^{1+\e}h^{-1},
$$
so that the bound in (3.8) becomes
$$
hX^{2\s(k)+2\e} + X^{1-\e}h^2 \ll X^{1-\e}h^2
$$
for
$$
h \;\ge\; X^{2\s(k)-1+3\e}.
$$
Therefore (2.4) follows if $0 < \e_1 < {1\over3} (\t(k) - 2\s(k)+1)$. This completes
the proof of Theorem 1.

\medskip
\section{Proof of Theorem 2}
\medskip
We may suppose that $X\ge h\;(\ge2)$ are integers, for otherwise note that by replacing
$X$ with $[X]$ in (2.5) and (2.6) we make an error which is, by trivial estimation,
$\ll_\e h^2X^\e$, and likewise for $h$. Write
$$
I_k(X,h) = \int_X^{2X}\Bigl(\D_k(x+h) - \D_k(x)\Bigr)^2\d x
= \sum_{X\le m\le 2X-1}I_{k,h}(m),\leqno(4.1)
$$
say, where for $m\in\NN$ we set
$$
I_{k,h}(m) \;:=\; \int_m^{m+1-0}\Bigl(\D_k(x+h) - \D_k(x)\Bigr)^2\d x.
$$
Recall that
$$
\sum_{n\le y}d_k(n) = \sum_{n\le [y]}d_k(n)\qquad(y>1),
$$
so that
$$
\sum_{x<n\le x+h}d_k(n) = \sum_{m<n\le m+h}d_k(n)\qquad(m \le x < m+1;\,m,h\in\NN).
$$
If $Q_{k-1} := P_{k-1}+P'_{k-1}$ (see (1.1)), then we have (for some $0\le \t \le1$)
$$
\leqno(4.2)
$$
\begin{eqnarray*}
&{}&\D_k(x+h) - \D_k(x) = \sum_{x<n\le x+h}d_k(n)\\ &+& xP_{k-1}(\log x)
- (x+h)P_{k-1}(\log(x+h))\\&
=& \sum_{m<n\le m+h}d_k(n) - hQ_{k-1}(\log(x + \t h))\\&
=& \sum_{m<n\le m+h}d_k(n) - hQ_{k-1}(\log m) + O(h^2X^{-1}\log^{k-1}X)\\&
= &\D_k(m+h) - \D_k(m)+ O(h^2X^{-1}\log^{k-1}X),
\end{eqnarray*}
where we used the mean value theorem. By using (4.2) we see that the left-hand side of
(4.1) becomes
$$
\leqno(4.3)
$$
\begin{eqnarray*}
&{}&
\sum_{X\le m\le 2X-1}I_{k,h}(m) \\&=& \sum_{X\le m\le 2X-1}
\int\limits_m^{m+1-0}{\left(\D_k(m+h) - \D_k(m)+
 O(h^2X^{-1}\log^{k-1}X)\right)}^2\d x \\&
= & \sum_{X\le m\le 2X}
{\left(\D_k(m+h) - \D_k(m)\right)}^2 + O_\e(X^\e h^2) \;\\&+&
  O\Bigl(\sum_{X\le m\le 2X-1}
\left|\D_k(m+h) - \D_k(m)\right|h^2X^{-1}\log^{k-1}X\Bigr)\\
&+& O(h^4X^{-1}\log^{2k-2}X).
\end{eqnarray*}
If $k\ge3$, we use the trivial bound (1.5) to obtain
$$
\sum_{X\le m\le 2X-1}
\left|\D_k(m+h) - \D_k(m)\right|h^2X^{-1}\log^{k-1}X \ll_\e h^3X^\e,
$$
and (2.5) follows, since
$$
h^4X^{-1}\log^{2k-2}X \ll h^2\log^{2k-2}X\quad(1 \ll h = h(X) \le \hf\sqrt{X}),
$$
and the last term is smaller than all the error terms in (2.5) and (2.6).

\medskip
If $k=2$, then we apply the Cauchy-Schwarz inequality to the sum in the
first $O$-term in (4.3) and use (1.7). We obtain
$$
\leqno(4.4)
$$
\begin{eqnarray*}
&{}&
\sum_{X\le m\le 2X-1}\left|\D(m+h) - \D(m)\right|h^2X^{-1}\log X \\&
\ll& X^{-1/2}h^2\log X\left(\sum_{X\le m\le 2X}(\D(m+h) - \D(m))^2\right)^{1/2}\\&
\ll& h^{5/2}\log^{5/2}X,
\end{eqnarray*}
and (2.5) follows from (4.3) and (4.4). This ends the proof of Theorem 2.

\pagebreak

\end{document}